\documentclass[12pt,leqno,twoside]{article}

\usepackage{amsmath,amssymb,amsthm}
\usepackage{titlesec,hyperref}
\usepackage{fancyhdr}

\usepackage{color}

\usepackage[margin=3.17cm]{geometry} 

\usepackage{fancyhdr}

\titleformat{\subsection}{\it}{\thesubsection.\enspace}{1pt}{}

\newtheorem{theo}{Theorem}[section]
\newtheorem{lemm}[theo]{Lemma}
\newtheorem{defi}[theo]{Definition}

\newtheorem{prop}[theo]{Proposition}
\newtheorem{rema}[theo]{Remark}
\numberwithin{equation}{section}

\def\ep{\varepsilon}

\allowdisplaybreaks
\begin{document}
\title{Global existence and well-posedness of 2D viscous shallow water system in Sobolev spaces with low regularity
\hspace{-4mm}
}

\author{Yanan Liu$^1$
\quad Zhaoyang Yin$^2$ \\[10pt]
Department of Mathematics, Sun Yat-sen University,\\
510275, Guangzhou, P. R. China\\[5pt]
}
\footnotetext[1]{Email: \it babyinarm@126.com}
\footnotetext[2]{Email: \it mcsyzy@mail.sysu.com.cn}
\date{}
\maketitle

\begin{abstract}
In this paper we consider the Cauchy problem for 2D viscous
shallow water system in $H^s(\mathbb{R}^2)$, $s>1$. We first prove
the local well-posedness of this problem by using the
Littlewood-Paley theory, the Bony decomposition, and the theories
of transport equations and transport diffusion equations. Then, we
get the global existence of the system with small initial data in
$H^s(\mathbb{R}^2)$, $s>1$. Our obtained result improves the
recent result in \cite{W}.

\vspace*{5pt}

\vspace*{5pt} \noindent{\it Keywords}: Viscous shallow water
system; the Littlewood-Paley theory; the Bony Decomposition;
Sobolev spaces; global existence
\end{abstract}

\vspace*{10pt}

\tableofcontents

\section{Introduction}
~~We consider the following Cauchy problem for the following 2D
viscous shallow water equations
\begin{align}
\left\{
\begin{array}{l}
h(u_t+(u\cdot\nabla)u)-\nu\nabla\cdot(h\nabla u)+h\nabla h=0, \\[1ex]
h_t+div(hu)=0,  \\[1ex]
u|_{t=0}=u_0,\quad h|_{t=0}=h_0,
\end{array}
\right.
\end{align}
where $h(x,t)$ is the height of fluid surface,
$u(x,t)=(u^1(x,t),u^2(x,t))$ is the horizontal velocity field,
$x=(x_1,x_2)\in{\mathbb{R}}^2,$ and $0<\nu<1$ is the viscous
coefficient. For the initial data $h_0(x)$, we suppose that it is
a small perturbation of some positive constant $\bar h_0$. In the
paper, we mainly study the Cauchy problem (1.1) in Sobolev spaces
$H^s({\mathbb{R}}^2)$, $s>1$. For the sake of convenience, we let
the notation $B^s_{p,r}$ stand for $B^s_{p,r}({\mathbb{R}}^2)$ in
the following text, and also let the notations $L^p$ and $H^s$
stand for $L^p({\mathbb{R}}^2)$ and $H^s({\mathbb{R}}^2)$,
respectively.

Recently, Bresch et al. \cite{D, D.B} have systematically
introduced the viscous shallow water equations. Bui in \cite{bui}
proved the local existence and uniqueness of classical solutions
to the Cauchy-Dirichlet problem for the shallow water equations
with initial data in $C^{2+\alpha}$ by using Lagrangian
coordinates and H\"{o}lder space estimates. Kloeden in \cite{K}
and Sundbye in \cite{S1} independently showed the global existence
and uniqueness of classical solutions to the Cauchy-Dirichlet
problem using Sobolev space estimates by following the energy
method of Matsumura and Nishida \cite{M-N}. Sundbye in \cite{S2}
proved the existence and uniqueness of classical solutions to the
Cauchy problem using the method of \cite{M-N}. Wang and Xu in in
\cite{W}, obtained local solutions for any initial data and global
solutions for small initial data $h_0-\bar h_0,u_0\in H^s$, $s>2$.
Haspot got global existence in time for small initial data $h_0$,
$h_0-\bar h_0\in \dot{B}^{0}_{2,1}\cap \dot{B}^{1}_{2,1}$ and $u_0
\in \dot{B}^{0}_{2,1}$ as a special case in \cite{H}, and Chen,
Miao and Zhang in \cite{C-M-Z} to prove the local existence in
time for general initial data and the global existence in time for
small initial data where $h_0-\bar h_0\in \dot{B}^{0}_{2,1}\cap
\dot{B}^{1}_{2,1}$ and $u_0 \in \dot{B}^{0}_{2,1}$ with additional
conditions that $h\geq h_0$ and $h_0$ is a strictly positive
constant.

In the paper, we mainly use the Littlewood-Paley theory, the Bony
decomposition and the Besov space theories for transport equations
and transport-diffusion equations to obtain the local existence
and uniqueness of solutions for any initial data in $H^s,\,s>1$.
And thus we can get the global existence of the system (1.1) with
the small enough initial data in $H^s,\,s>1$. The main result of
this paper is as follows.
\begin{theo}\label{t1}
Let $u_0,~h_0-\bar{h}_0\in H^s$, $s>1$,
$\|h_0-\bar{h}_0\|_{H^s}<<\bar{h}_0$. Then there exists a positive
time $T$, a unique solution $(u,h)$ of the Cauchy problem (1.1)
such that
$$u,h-\bar{h}_0\in L^{\infty}([0,T];H^s),~u\in L^2([0,T];H^{s+1}).$$
Moreover, there exists a constant $\eta$, such that if $\|u_0\|_{H^s}+\|h_0\|_{H^s}\leq \eta$, then the time $T$ can be infinity.
\end{theo}

\section{Preliminaries}

First of all, we transform the system (1.1). For a sake of convenience, we take $\bar{h}_0=1$. Substituting $h$ by $1+h$ in (1.1), we have
\begin{align}
\left\{
\begin{array}{l}
u_t+(u\cdot\nabla)u-\nu\Delta u-\nu\nabla(\ln(1+h))\nabla u+\nabla h=0, \\[1ex]
h_t+div\,u+div(hu)=0,  \\[1ex]
u|_{t=0}=u_0,\quad h|_{t=0}=h_0,
\end{array}
\right.
\end{align}
~~~~~~~~~~~~~~~~~~~~~~~~~~~$x\in \mathbb{R}^2$, $t\in[0,\infty).$

here $h_0\in H^s$, $\|h_0\|_{H^s}\leq \frac{1}{8C_0}$, and $C_0$
will be determined below.

Next we introduce the well-known Littlewood-Paley decomposition
briefly.

\begin{prop}\label{l2} \cite{H.J}
Littlewood-Paley decomposition:

\noindent Let $\mathcal{B}=\{\xi\in{\mathbb{R}^2},~
|\xi|\leq\frac{4}{3}\}$ be a ball, and let
$\mathcal{C}=\{\xi\in{\mathbb{R}^2},
~\frac{3}{4}\leq|\xi|\leq\frac{8}{3} \}$ be an annulus. There
exist two radial functions $\chi$ and $\varphi$ valued in the
interval $[0,1]$, belonging respectively to
$\mathcal{D}(\mathcal{B})$ and $\mathcal{D}(\mathcal{C})$, such
that
\begin{align}
\forall~\xi\in{\mathbb{R}}^2,~~\chi(\xi)+\Sigma_{j\geq0}\varphi(2^{-j}\xi)=1,
\end{align}
\begin{align}
\forall~\xi\in{\mathbb{R}}^2\backslash\{0\},~~\Sigma_{j\in\mathbb{Z}}\varphi(2^{-j}\xi)=1,
\end{align}
\begin{align}
|j-j'|\geq2~\Rightarrow Supp\varphi(2^j\cdot)\cap Supp\varphi(2^{j'}\cdot)=\emptyset,
\end{align}
\begin{align}
j\geq2~\Rightarrow Supp\chi\cap Supp\varphi(2^j\cdot)=\emptyset,
\end{align}

the set $\tilde{\mathcal{C}}\overset{def}{=}~B(0,2/3)+\mathcal{C}$ is an annulus, and we have
\begin{align}
|j-j'|\geq5~\Rightarrow 2^j\tilde{\mathcal{C}}\cap 2^{j'}\mathcal{C}=\emptyset.
\end{align}

Further, we have
\begin{align}
\forall~\xi\in{\mathbb{R}}^2,~~\frac{1}{2}\leq\chi^2(\xi)+\Sigma_{j\geq0}\varphi^2(2^{-j}\xi)\leq1,
\end{align}
\begin{align}
\forall~\xi\in{\mathbb{R}}^2\backslash\{0\},~~\frac{1}{2}\leq\Sigma_{j\in\mathbb{Z}}\varphi^2(2^{-j}\xi)\leq1.
\end{align}
\end{prop}

Now we can define the nonhomogeneous dyadic blocks $\Delta_j$ and the nonhomogeneous low-frequency cut-off operator $S_j$ as follows:
$$\Delta_j u=0,~if~j\leq-2,~~\Delta_{-1}u=\chi(D)u=\int_{\mathbb{R}^2}\tilde{h}(y)u(x-y)dy,$$
$$\Delta_ju=\varphi(2^{-j}D)u=2^{jd}\int_{\mathbb{R}^2}h(2^jy)u(x-y)dy~~if~~j\geq0,$$
and
$$S_ju=\sum_{j'\leq j-1}\Delta_{j'}u,$$
where $h=\mathcal{F}^{-1}\varphi$ and
$\tilde{h}=\mathcal{F}^{-1}\chi$.

Next we define the Besov spaces:
\begin{defi} \label{l3} \cite{H.J}
Let $s\in \mathbb{R}$ and $(p,r)\in[1,\infty]^2$. The
nonhomogeneous Besov space $B^s_{p,r}$ consists of all tempered
distribution $u$ such that:
$$\left(\sum_{j\geq-1}(2^{js}\|\Delta_ju\|_{L^p})\right)_{\ell^r}<\infty,$$
and naturally the Besov norm is defined as follows
$$\|u\|_{B^s_{p,r}}=\left(\sum_{j\geq-1}(2^{js}\|\Delta_ju\|_{L^p})\right)_{\ell^r}.$$
Particularly, when $p,r=2$, the Besov space is coincident with the
Sobolev space, i.e:
$$B^s_{2,2}=H^s.$$
\end{defi}
\begin{defi}\label{l4} \cite{H.J}
The Bony decomposition: The nonhomogeneous paraproduct of $v$ by
$u$ is defined by
$$T_uv=\underset{j}{\sum}S_{j-1}u\Delta_jv.$$
The nonhomogeneous remainder of $u$ by $v$ is defined by
$$R(u,v)=\underset{|k-j|\leq1}{\sum}\Delta_ku\Delta_jv.$$
The operators $T$ and $R$ are bilinear, and we have the following Bony decomposition
$$uv=T_vu+T_uv+R(u,v).$$
\end{defi}
\begin{lemm}\label{l5} \cite{H.J}
Bernstein-Type inequalities:

Let $\mathcal{C}$ be an annulus and $\mathcal{B}$ a ball. A constant $C$ exists such that for any nonnegative integer $k$, any couple $(p,q)$ in $[1,\infty]^2$ with $q\geq p\geq1$, and any function u of $L^p$, we have
$$Supp\,\hat{u}\subset \lambda\mathcal{B}\Rightarrow\,\underset{|\alpha|=k}{sup}\,\|\partial^{\alpha}u\|_{L^q}
\leq C^{k+1}\lambda^{k+d(\frac{1}{p}-\frac{1}{q})}\|u\|_{L^p},$$
$$Supp\,\hat{u}\subset \lambda\mathcal{C}\Rightarrow\,C^{-k-1}\lambda^k\|u\|_{L^p}\leq\underset{|\alpha|=k}{sup}\,\|\partial^{\alpha}u\|_{L^q}
\leq C^{k+1}\lambda^k\|u\|_{L^p}.$$
\end{lemm}
Then we give some properties of the Besov spaces which will be used in this paper.
\begin{lemm}\label{l6} \cite{H.J}
Let $1\leq p_1\leq p_2\leq\infty$ and $1\leq r_1\leq
r_2\leq\infty$. Then for any $s\in \mathbb{R}$, the space
$B^s_{p_1,r_1}$ is continuously embedded in
$B^{s-d(\frac{1}{p_1}-\frac{1}{p_2})}_{p_2,r_2}$. Obviously, we
also have that the space $B^{s_2}_{p,r}$ is continuously embedded
in $B^{s_1}_{p,r}$ and $B^{s_2}_{p,\infty}$ is continuously
embedded in $B^{s_1}_{p,1}$ if $s_1<s_2$.
\end{lemm}
\begin{lemm}\label{0} \cite{H.J}
If $s_1$ and $s_2$ are real numbers such that $s_1<s_2$,
$\theta\in(0,1)$, and $p,r\in[1,\infty]$, then we have
$$\|u\|_{B^{\theta s_1+(1-\theta)s_2}_{p,r}}\leq\|u\|^\theta_{B^{s_1}_{p,r}}\|u\|^{1-\theta}_{B^{s_2}_{p,r}}.$$
\end{lemm}
\begin{lemm}\label{l7} \cite{H.J}
If $u\in B^s_{p,r}$, then $\nabla u\in B^{s-1}_{p,r}$, and we have
$$\|\nabla u\|_{B^{s-1}_{p,r}}\leq C\|u\|_{B^s_{p,r}}.$$
\end{lemm}
\begin{lemm}\label{l8} \cite{H.J}
The set $B^s_{p,r}$ is a Banach space. and satisfies the Fatou
property, namely, if $(u_n)_{n\in N}$ is a bounded sequence of
$B^s_{p,r}$. Then an element $u$ of $B^s_{p,r}$ and a subsequence
$u_{\psi(n)}$ exist such that:

$\underset{n\rightarrow\infty}{lim}~u_{\psi(n)}=u$ in $\mathcal{S}'$,  $\|u\|_{B^s_{p,r}}\leq C \underset{n\rightarrow\infty}{lim} inf\|u_{\psi(n)}\|_{B^s_{p,r}}$.
\end{lemm}
\begin{lemm}\label{l9} \cite{H.J}
If $s>\frac{d}{p}$ or $s=\frac{d}{p},~r=1$, then the $B^s_{p,r}$
space is continuously embedded in $L^{\infty},$ i.e
$$\|u\|_{L^{\infty}}\leq C_{s,p}\|u\|_{B^s_{p,r}}.$$
\end{lemm}

\begin{lemm}\label{l10} \cite{H.J}
Let f be a smooth function, $f(0)=0$, $s>0,~(p,r)\in[1,\infty]^2$.
If $u\in B^s_{p,r}\cap L^{\infty}$, then so does $f\circ u$, and
we have
$$\|f\circ u\|_{B^s_{p,r}}\leq C\left(s,f',\|u\|_{L^{\infty}}\right)\|u\|_{B^s_{p,r}}.$$
\end{lemm}
\begin{lemm}\label{l11} \cite{H.J}
A constant $C$ exists which satisfies the following inequalities
for any couple of real numbers $(s,t)$ with t negative and any
$(p,r_1,r_2)$ in $[1,\infty]^3$:
$$\|T\|_{\mathcal{L}(L^{\infty}\times B^s_{p,r};B^s_{p,r})}\leq C^{|s|+1},$$
$$\|T\|_{\mathcal{L}(B^t_{\infty,r_1}\times B^s_{p,r_2};B^{s+t}_{p,r})}\leq \frac{C^{|s+t|+1}}{-t}
~~with~~\frac{1}{r}\overset{def}{=}min\{1,\frac{1}{r_1}+\frac{1}{r_2}\}.$$
\end{lemm}
\begin{lemm}\label{l12} \cite{H.J}
A constant $C$ exists which satisfies the following inequalities.
Let $(s_1,s_2)$ be in $\mathbb{R}^2$ and $(p_1,p_2,r_1,r_2)$ be in
$[1,\infty]^4$. Assume that
$$\frac{1}{p}\overset{def}{=}\frac{1}{p_1}+\frac{1}{p_2}\leq1~~and~~
\frac{1}{r}\overset{def}{=}\frac{1}{r_1}+\frac{1}{r_2}\leq1.$$
If $s_1+s_2>0$, then we have, for any $(u,v)$ in $B^{s_1}_{p_1,r_1}\times B^{s_2}_{p_2,r_2}$,
$$\|R(u,v)\|_{B^{s_1+s_2}_{p,r}}\leq\frac{C^{|s_1+s_2|+1}}{s_1+s_2}\|u\|_{B^{s_1}_{p_1,r_1}}\|v\|_{B^{s_2}_{p_2,r_2}}.$$
If $r=1$ and $s_1+s_2=0$, then we have, for any $(u,v)$ in $B^{s_1}_{p_1,r_1}\times B^{s_2}_{p_2,r_2}$,
$$\|R(u,v)\|_{B^0_{p,\infty}}\leq C^{|s_1+s_2|+1}\|u\|_{B^{s_1}_{p_1,r_1}}\|v\|_{B^{s_2}_{p_2,r_2}}.$$
\end{lemm}

\begin{lemm}\label{l13} \cite{H.J}
For any $s>0$ and $(p,r)\in[1,\infty]^2$, the space $B^s_{p,r}\cap
L^{\infty}$ is an algebra, and a constant exists such that:

$$\|uv\|_{B^s_{p,r}}\leq \frac{C^{s+1}}{s}\Big(\|u\|_{L^{\infty}}\|v\|_{B^s_{p,r}}+\|v\|_{L^{\infty}}\|u\|_{B^s_{p,r}}\Big).$$
Moreover, if $s>\frac{d}{p}~or~s=\frac{d}{p},r=1$, we have
$$\|uv\|_{B^s_{p,r}}\leq \frac{C^{s+1}}{s}\|u\|_{B^s_{p,r}}\|v\|_{B^s_{p,r}}.$$
\end{lemm}

For the transport equations
\begin{align}
\left\{
\begin{array}{l}
\partial_tf+v\cdot\nabla f=g,\\
f_{|t=0}=f_0,
\end{array}
\right.
\end{align}
we have
\begin{lemm}\label{l14} \cite{H.J}
Let $1\leq p\leq p_1\leq\infty,~1\leq r\leq\infty$. Assume that
\begin{align}
s\geq-d\,min\left(\frac{1}{p_1},\frac{1}{p'}\right) \quad or \quad s\geq-1-d\,min\left(\frac{1}{p_1},\frac{1}{p'}\right)~if~div\,v=0
\end{align}
with strict inequality if $r<\infty$.

There exists a constant $C$, depending only on $d, p, p_1, r$ and
$s$, such that for all solutions $f\in
L^{\infty}([0,T];B^s_{p,r})$ of (2.9), initial data $f_0$ in
$B^s_{p,r}$, and $g$ in $L^1([0,T];B^s_{p,r})$, we have, for a.e.
$t\in[0,T]$,
$$\|f\|_{\tilde{L}_t^{\infty}(B^s_{p,r})}\leq\left(\|f_0\|_{B^s_{p,r}}+
\int^t_0exp(-CV_{p_1}(t'))\|g(t')\|_{B^s_{p,r}}dt'\right)exp(CV_{p_1}(t))$$
with, if the inequality is strict in (2.10),
\begin{align}
V'_{p_1}(t)=\left\{\begin{array}{l}\|\nabla v(t)\|_{B^{s-1}_{p_1,r}},~if~s>1+\frac{d}{p_1}~or~s=1+\frac{d}{p_1},~r=1,\\
\|\nabla v(t)\|_{B^{\frac{d}{p_1}}_{p_1,\infty}\cap
L^{\infty}},~if~s<1+\frac{d}{p_1},
\end{array}\right.
\end{align}
and if the equality holds in (2.10) and $r=\infty$,
$$V'_{p_1}=\|\nabla v(t)\|_{B^{\frac{d}{p_1}}_{p_1,1}}.$$
If $f=v$, then for all $s>0$ ($s>-1,~if~div\,u=0$), the estimate
holds with
$$V'_{p_1}(t)=\|\nabla u\|_{L^{\infty}},$$
where $\|u\|_{\tilde{L}^{\rho}_T(B^s_{p,r})}$ is defined in Lemma
\ref{l15}.
\end{lemm}

For the transport diffusion equations
\begin{align}
\left\{
\begin{array}{l}
\partial_tf+v\cdot\nabla f-\nu\Delta f=g,\\
f_{|t=0}=f_0,
\end{array}
\right.
\end{align}
we have the following lemma.
\begin{lemm}\label{l15} \cite{H.J}
Let $1\leq p_1\leq p\leq\infty,~1\leq r\leq\infty,~s\in\mathbb{R}$
satisfy (2.10), and let $V_{p_1}$ be defined as in Lemma \ref{l8}.

There exists a constant $C$ which depends only on $d, r, s$ and $s-1-\frac{d}{p_1}$ and is such that for any smooth solution of (11) and $1\leq\rho_1\leq\rho\leq\infty,$ we have
$$\nu^{\frac{1}{\rho}}\|f\|_{\tilde{L}^{\rho}_T(B^{s+\frac{2}{\rho}}_{p,r})}\leq Ce^{C(1+\nu T)^{\frac{1}{\rho}}V_{p_1}(T)}\Big((1+\nu T)^{\frac{1}{\rho}}\|f_0\|_{B^s_{p,r}}$$
$$+(1+\nu T)^{1+\frac{1}{\rho}-\frac{1}{\rho_1}}\nu^{\frac{1}{\rho_1}-1}\|g\|_
{\tilde{L}^{\rho_1}_T(B^{s-2+\frac{2}{\rho_1}}_{p,r})}\Big).$$
\end{lemm}

For the space $\tilde{L}^{\rho}_T(B^s_{p,r})$, we have the following properties:
\begin{lemm}\label{l16} \cite{H.J}
For all $T>0,~s\in\mathbb{R},$ and $1\leq r,\rho\leq\infty$, we
set
$$\|u\|_{\tilde{L}^{\rho}_T(B^s_{p,r})}\overset{def}{=}\|2^{js}\|\Delta_ju\|_{L^{\rho}_T(L^p)}\|
_{l^r(\mathbb{Z})}.$$
We can then define the space $\tilde{L}^{\rho}_T(B^s_{p,r})$ as the set of tempered distributions
$u$ over $(0,T)\times \mathbb{R}^d$ such that $\|u\|_{\tilde{L}^{\rho}_T(B^s_{p,r})}\leq\infty$.
By the Minkowski inequality, we have
$$\|u\|_{\tilde{L}^{\rho}_T(B^s_{p,r})}\leq\|u\|_{L^{\rho}_T(B^s_{p,r})}\quad if~r\geq\rho,$$
$$\|u\|_{L^{\rho}_T(B^s_{p,r})}\leq\|u\|_{\tilde{L}^{\rho}_T(B^s_{p,r})}\quad if~r\leq\rho.$$
The general principle is that all the properties of continuity for the product, composition, remainder, and
paraproduct remain true in those spaces.

Moreover when $s>0,~1\leq p\leq\infty,~1\leq \rho,\rho_1,\rho_2,\rho_3,\rho_4\leq\infty,$ and
$$\frac{1}{\rho}=\frac{1}{\rho_1}+\frac{1}{\rho_2}=\frac{1}{\rho_3}+\frac{1}{\rho_4},$$
we have
$$\|uv\|_{\tilde{L}^{\rho}_T(B^s_{p,r})}\leq C\left(\|u\|_{\tilde{L}^{\rho_1}_T(L^{\infty})}
\|v\|_{\tilde{L}^{\rho_2}_T(B^s_{p,r})}+\|v\|_{\tilde{L}^{\rho_3}_T(L^{\infty})}
\|u\|_{\tilde{L}^{\rho_4}_T(B^s_{p,r})}\right).$$
\end{lemm}

\begin{rema}
 For the sake of convenience, for the fixed $s,p,r$, we let $C_0(\geq1)$ be the maximal constant of Lemmas 2.5-2.16.
\end{rema}
\section{The local existence}
First of all, for the sake of convenience, we only consider the case of $s<2$, the case of $s>2$ was proved by Wang and Xu in \cite{W}, the case of $s=2$ is similar to $s<2$ but more easier.
In order to study the local existence of solution, we define the function set $(u,h)\in\chi([0,T], s, E_1, E_2)$,  if

$$u\in\tilde{L}^{\infty}([0,T];H^s)\cap\tilde{L}^2([0,T];H^{s+1})\cap\tilde{L}^{\frac{4}{3-s}}([0,T];H^{\frac{s+3}{2}}),$$ $$h\in\tilde{L}^{\infty}([0,T];H^s),$$
and
$$\|u\|_{\tilde{L}^{\infty}([0,T];H^s)}\leq E_1,~\|u\|_{\tilde{L}^2([0,T];H^{s+1})}\leq E_1,~
\|u\|_{\tilde{L}^{\frac{4}{3-s}}([0,T];H^{\frac{s+3}{2}})}\leq E_1,$$
$$\|h\|_{\tilde{L}^{\infty}([0,T];H^s)}\leq E_2,$$
where
 $$E_1=8\nu^{-1}C_0\|u_0\|_{H^s},~E_2=4\|h_0\|_{H^s}.$$

Next, we will prove Theorem \ref{t1} by the method of successive approximations.
Let us define the sequence $(u_n,h_n)$ by the following linear system:
\begin{align}
\left\{
\begin{array}{l}
(u_1,h_1)=S_2(u_0,h_0),\\[1ex]
\partial_tu_{n+1}+(u_n\cdot\nabla)u_{n+1}-\nu\Delta u_n=\frac{\nu}{1+h_n}\nabla h_n\nabla u_n+\nabla h_n,\\[1ex]
\partial_th_{n+1}+(u_n\cdot\nabla)h_{n+1}=-div\,u_n-h\,div\,u_n,\\[1ex]
(u_{n+1},h_{n+1})_{|t=0}=S_{n+2}(u_0,h_0).
\end{array}
\right.
\end{align}

Since $S_q$ are smooth operators, the initial date $S_{n+2}(u_0,h_0)$ are smooth functions. If $(u_n,h_n)\in \chi ([0,T],s,E_1,E_2)$ are smooth, then we have that for any $t\in[0,T]$,
$$\|h_n\|_{L^\infty} \leq C_0 \|h_n\|_{H^s} \leq C_0E_2=4C_0\|h_0\|_{H^s} \leq \frac{4C_0}{8C_0}=\frac{1}{2}.$$
Thus $\frac{\nu}{1+h_n}\nabla h_n\nabla u_n+\nabla h_n$ and $-div\,u_n-h\,div\,u_n$ are also smooth functions.
Note that the first equation in (3.1) is a transport diffusion equation for $u_{n+1}$, and the second equation is a transport
equation for $h_{n+1}$. Then the local existence of the smooth function for the Cauchy problem (3.1) is obvious.

 We split the proof of Theorem 1.1 into two steps: Estimation for big norms and Convergence for small norms.
\subsection{Estimation for big norm}
In this subsection, we want to prove the following proposition.
\begin{prop}\label{p1}
Suppose that $(u_0,h_0)\in H^s\times H^s$, $s\in(1,2)$ and $\|h_0\|_{B^s_{p,r}}\leq \frac{1}{8C_0}$, then
 there exists a positive time $T$, such that for any $n\in N$, $(u_n,h_n)\in\chi ([0,T_1],s,E_1,E_2)$.
\end{prop}
Proof: Let $T\in(0,1]$ satisfy firstly
$$e^{2C_0^2E_1T}\leq2,~(1+\nu T)\leq2,~T\leq \frac{9}{64}.$$
 Then we prove the proposition by induction. Firstly letting $(u_1,h_1)=S_2(u_0,h_0)$, thus we have
$$\|u_1\|_{\tilde{L}^{\infty}_T(H^s)}\leq\|u_0\|_{H^s}\leq E_1,~
\|h_1\|_{\tilde{L}^{\infty}_T(H^s)}\leq\|h_0\|_{H^s}\leq E_2.$$
By Proposition \ref{l2}, we also have
$$\|u_1\|_{L^2_T(H^{s+1})}\leq T^{\frac{1}{2}}\|S_2u_0\|_{L^\infty(H^{s+1})}\leq \frac{32}{3}T^{\frac{1}{2}}\|S_2u_0\|_{H^s}\leq E_1.$$

We suppose now that $$\|u_n\|_{\tilde{L}^{\infty}_T(H^s)}\leq E_1,~\|u_1\|_{L^2_T(H^{s+1})}\leq E_1 ,~
\|h_n\|_{\tilde{L}^{\infty}_T(H^s)}\leq E_2.$$

Then for $u_{n+1}$, in the view of Lemma \ref{l15}, taking a positive real number $\rho$ such that $\frac{2}{\rho}<s-1$, then for all $t\leq T$, we have
\begin{align}
\begin{array}{l}
\|u_{n+1}\|_{\tilde{L}_t^{\infty}(H^s)}\leq C_0e^{C_0\int_0^t\|\nabla u_n(t')\|_{H^1\cap L^\infty}dt'}\Big(\|S_{n+2}u_0\|_{H^s}+(1+\nu T)^{1-\frac{1}{\rho}}\nu^{\frac{1}{\rho}-1}\times\\[1ex]
\|\nu\nabla(\ln(1+h_n))\nabla u_n-\nabla h_n\|_{\tilde{L}^{\rho}_T(H^{s-2+\frac{2}{\rho}})}\Big)
\\[1ex]
\leq2C_0\big(\frac{E_1}{4C_0}+2\nu^{\frac{1}{\rho}-1}\|\nabla h_n\|_{\tilde{L}_t^\rho(H^{s-2+\frac{2}{\rho}})}+
2\nu^{\frac{1}{\rho}}\|\nabla(\ln(1+h_n))\nabla u_n\|_{\tilde{L}_t^{\rho}(H^{s-2+\frac{2}{\rho}})}\big)\\[1ex]
\leq
\frac{E_1}{2}+CT^{\frac{1}{\rho}}\|h_n\|_{\tilde{L}^\infty_T(H^s)}+C\|\nabla(\ln(1+h_n))\nabla
u_n\|_{\tilde{L}_t^{\rho}(H^{s-2+\frac{2}{\rho}})}.\end{array}
\end{align}
Now, we estimate $\|\nabla(\ln(1+h_n))\nabla u_n\|_{\tilde{L}_t^{\rho}(H^{s-2+\frac{2}{\rho}})}$. By Lemmas \ref{l11} and \ref{l12}, we get
\begin{align}
\begin{array}{l}
\|\nabla(\ln(1+h_n))\nabla u_n\|_{\tilde{L}_t^{\rho}(H^{s-2+\frac{2}{\rho}})}\\[1ex]
\leq\|T_{\nabla(\ln(1+h_n))}\nabla u_n+T_{\nabla u_n}\nabla(\ln(1+h_n))+R\big(\nabla(\ln(1+h_n)),\nabla u_n\big)\|_{L_t^{\rho}(H^{s-2+\frac{2}{\rho}+\ep})}\\[1ex]
\leq CT^{\frac{1}{\rho}}\big(\|\nabla(\ln(1+h_n))\|_{L^{\infty}_T(B^{s-2}_{\infty,\infty})}\|\nabla u_n\|
_{L^{\infty}_T(H^{\frac{2}{\rho}+\ep})}+\|\nabla u_n\|_{L^{\infty}_T(B^{s-2}_{\infty,\infty})}
\|\nabla(\ln(1+h_n))\|_{L^{\infty}_T(H^{\frac{2}{\rho}+\ep})}\big)\\[1ex]
+CT^{\frac{1}{\rho}}\|R\big(\nabla(\ln(1+h_n)),\nabla u_n\big)\|_{L^{\infty}_T(B^{s-1+\frac{2}{\rho}+\ep}_{1,2})}\\[1ex]
\leq CT^{\frac{1}{\rho}}\|u_n\|_{L^{\infty}_T(H^s)}\|h_n\|_{L^{\infty}_T(H^s)}+CT^{\frac{1}{\rho}}
\|\nabla(\ln(1+h_n))\|_{L^{\infty}_T(H^{s-1})}\|\nabla u_n\|_{L^{\infty}_T(H^{\frac{2}{\rho}+\ep})}\\[1ex]
\leq CT^{\frac{1}{\rho}}\|u_n\|_{L^{\infty}_T(H^s)}\|h_n\|_{L^{\infty}_T(H^s)},
\end{array}
\end{align}
where $\ep$ is a small enough positive real number, and it does't affect the direction of the inequality:
$$if~ 1+\frac{2}{\rho}<s,\, then~ 1+\frac{2}{\rho}+\ep<s.$$
Thus combining with (3.2), we get
\begin{align}
\begin{array}{l}
\|u_{n+1}\|_{\tilde{L}_t^{\infty}(H^s)}\leq \frac{E_1}{2}+CT^{\frac{1}{\rho}}\|h_n\|_{\tilde{L}^\infty_T(H^s)}+
CT^{\frac{1}{\rho}}\|u_n\|_{L^{\infty}_T(H^s)}\|h_n\|_{L^{\infty}_T(H^s)}\\[1ex]
\leq \frac{E_1}{2}+CT^{\frac{1}{\rho}}E_2+CT^{\frac{1}{\rho}}E_1E_2.
\end{array}
\end{align}
\begin{rema}
In the calculation of (3.3), if $s=2$, we can estimate
$$\|T_{\nabla(\ln(1+h_n))}\nabla u_n\|_{L_t^{\rho}(H^{s-2+\frac{2}{\rho}+\ep})}$$
as follows:
\begin{align}
\begin{array}{l}
\|T_{\nabla(\ln(1+h_n))}\nabla u_n\|_{L_t^{\rho}(H^{s-2+\frac{2}{\rho}+\ep})}\\[1ex]
\leq CT^{\frac{1}{\rho}}\|\nabla(\ln(1+h_n))\|_{L^\infty_T(B^{s-2-\ep}_{\infty,\infty})}
\|\nabla u_n\|_{L^\infty_T(H^{\frac{2}{\rho}+2\ep})}\\[1ex]
\leq CT^{\frac{1}{\rho}}\|h_n\|_{L^\infty_T(H^s)}\|u_n\|_{L^\infty_T(H^s)}.
\end{array}
\end{align}
\end{rema}
Using the same method,  we deal with the similar case below as well.

For $h_{n+1}$, from Lemmas \ref{l14}, we obtain
\begin{align}
\begin{array}{l}
\|h_{n+1}\|_{\tilde{L}_t^{\infty}(H^s)}\leq e^{C_0\int_0^t\|\nabla u_n\|_{H^1\cap L^\infty}dt'}
\big(\|S_{n+2}h_0\|_{H^s}+\|div\,u_n\|_{\tilde{L}^1_T(H^s)}+\|h_n\,div\,u_n\|_{\tilde{L}^1_T(H^s)}\big)\\[1ex]
\leq2\|h_0\|_{H^s}+C\|u_n\|_{\tilde{L}^1_T(H^{s+1})}+C\|h_n\|_{L^{\infty}_T(H^s)}
\|div\,u_n\|_{L^1_T(H^s)}\\[1ex]
\leq \frac{E_2}{2}+(C+C\|h_n\|_{L^{\infty}_T(H^s)})T^{\frac{1}{2}}\|u_n\|_{L^2_T(H^{s+1})}\\[1ex]
\leq \frac{E_2}{2}+CT^{\frac{1}{2}}(1+E_2)E_1.
\end{array}
\end{align}
Then by Lemma \ref{l15} again, we have
\begin{align}
\begin{array}{l}
\|u_{n+1}\|_{L^2_t(H^{s+1})}\leq \nu^{-\frac{1}{2}}C_0e^{C_0(1+\nu T)^{\frac{1}{2}}\int_0^T\|\nabla u_n\|_{H^1\cap L^\infty}dt}\Big((1+\nu T)^{\frac{1}{2}}\|S_{n+1}u_0\|_{H^s}+\\[1ex]
(1+\nu T)\nu^{-\frac{1}{2}}\|\nu \nabla(\ln(1+h_n))\nabla u_n-\nabla h_n\|_{L^2_T(H^{s-1})}\Big)\\[1ex]
\leq 2C_0\nu^{-\frac{1}{2}}\big(2\|u_0\|_{H^s}+2\nu^{-\frac{1}{2}}C_0\|h_n\|_{L^2_T(H^s)}+2\nu^{-\frac{1}{2}}
\|\nabla(\ln(1+h_n))\nabla u_n\|_{L^2_T(H^{s-1})}\big)\\[1ex]
\leq 4C_0\nu^{-\frac{1}{2}}\|u_0\|_{H^s}+CT^{\frac{1}{2}}\|h_n\|_{L^\infty_T(H^s)}+C\|\nabla(\ln(1+h_n))\nabla u_n\|_{L^2_T(H^{s-1})}.
\end{array}
\end{align}
By Lemmas \ref{l11}-\ref{l12}, we have
\begin{align}
\begin{array}{l}
\|\nabla(\ln(1+h_n))\nabla u_n\|_{L^2_T(H^{s-1})}
\leq C\|\nabla(\ln(1+h_n))\|_{L^{\infty}_T(B^{s-2}_{\infty,\infty})}\|\nabla u_n\|_{L^2_T(H^1)}+\\[1ex]
C\|\nabla u_n\|_{L^2_T(L^\infty)}\|\nabla(\ln(1+h_n))\|_{L^{\infty}_T(H^{s-1})}+
C\|\nabla(\ln(1+h_n))\|_{L^{\infty}_T(B^{s-2}_{\infty,\infty})}\|\nabla u_n\|_{L^2_T(H^1)}\\[1ex]
\leq C\|h_n\|_{L^{\infty}_T(H^s)}\|u_n\|_{L^2_T(H^{\frac{s+3}{2}})}\\[1ex]
\leq CT^{\frac{s-1}{4}}\|h_n\|_{L^{\infty}_T(H^s)}\|u_n\|_{L^{\frac{4}{3-s}}_T(H^{\frac{s+3}{2}})}.
\end{array}
\end{align}
By Lemma \ref{0}, and taking $\theta=\frac{s-1}{2}$, we have
\begin{align}
\begin{array}{l}
\|u_n\|_{L^{\frac{4}{3-s}}_T(H^{\frac{s+3}{2}})}
\leq \|u_n\|^{\frac{s-1}{2}}_{L^\infty_T(H^s)}\|u_n\|
^{\frac{3-s}{2}}_{L^2_T(H^{s+1})}\leq E_1.
\end{array}
\end{align}
Combining (3.7), (3.8) and (3.9), we can obtain that
\begin{align}
\begin{array}{l}
\|u_{n+1}\|_{L^2_t(H^{s+1})}
\leq 4C_0\nu^{-\frac{1}{2}}\|u_0\|_{H^s}+CT^{\frac{1}{2}}\|h_n\|_{L^\infty_T(H^s)}\\[1ex]
+CT^{\frac{s-1}{4}}\|h_n\|_{L^{\infty}_T(H^s)}\|u_n\|^{\frac{s+1}{2}}_{L^\infty_T(H^s)}\|u_n\|
^{\frac{3-s}{2}}_{L^2_T(H^{s+1})}\\[1ex]
\leq \frac{E_1}{2}+CT^{\frac{1}{2}}E_2+CT^{\frac{s-1}{4}}E_2E_1,
\end{array}
\end{align}
where $C$ in (3.2)-(3.10) only depends on $C_0,~\nu$. Thus we only take a suitable $T$, such that
\begin{align}
\left\{
\begin{array}{l}
CT^{\frac{1}{\rho}}E_2\leq\frac{E_1}{4},~CT^{\frac{1}{\rho}}E_2\leq\frac{1}{4},\\
CT^{\frac{1}{2}}(1+E_2)E_1\leq \frac{E_2}{2},\\
CT^{\frac{1}{2}}E_2\leq\frac{E_1}{4},~CT^{\frac{s-1}{4}}E_2\leq\frac{1}{4}.
\end{array}
\right.
\end{align}
This completes the proof of Proposition \ref{p1}.
\subsection{Convergence of small norm}
\begin{prop}\label{p2}
Suppose that $(u_0,h_0)\in H^s\times H^s$, $s\in(1,2)$ , $\|h_0\|_{B^s_{p,r}}\leq \frac{1}{8C0}$, then
 there exists a positive time $T_1(\leq T)$, such that $(u_n,h_n)$ is a Cauchy sequence in $\chi ([0,T_1],s-1,E_1,E_2)$.
\end{prop}
Proof: From the equations in (3.1), we have
\begin{align}
\left\{
\begin{array}{l}
\partial_t(u_{n+1}-u_n)+(u_n\cdot\nabla)(u_{n+1}-u_n)-\nu\Delta(u_{n+1}-u_n)=\sum_{j=1}^5F_j\\[1ex]
\partial_t(h_{n+1}-h_n)+(u_n\cdot\nabla)(h_{n+1}-h_n)=\sum_{j=1}^4J_j,\\[1ex]
(u_{n+1}-u_n,h_{n+1}-h_n)_{|t=0}=\Delta_{n+1}(u_0,h_0),
\end{array}
\right.
\end{align}

where
\begin{align}
\begin{array}{l}
\sum_{j=1}^5F_j=(u_n-u_{n-1})\cdot\nabla u_n+\nabla(h_n-h_{n-1})+\frac{\nu}{1+h_n}\nabla h_n\nabla(u_n-u_{n-1})\\[1ex]
+\frac{\nu}{1+h_n}\nabla u_{n-1}\nabla(h_n-h_{n-1})+\nu
(\frac{1}{1+h_n}-\frac{1}{1+h_{n-1}})\nabla h_{n-1}\nabla u_{n-1},\\[2ex]
\sum_{j=1}^4J_j=(u_n-u_{n-1})\cdot\nabla h_n+div\,(u_n-u_{n-1})+h_n\,div\,(u_n-u_{n-1})\\[1ex]
+(h_n-h_{n-1})\,div\,u_{n-1}.
\end{array}
\end{align}

Then we estimate the Sobolev norms of $u_{n+1}-u_n$ and
$h_{n+1}-h_n$. For any $t\leq T_1\leq T$, by Lemma \ref{l15}, we
have
\begin{align}
\begin{array}{l}
\|u_{n+1}-u_n\|_{\tilde{L}^\infty_t(H^{s-1})}\leq C_0e^{\int_0^t\|\nabla u_n\|_{H^1\cap L^\infty}dt'}\times\\[1ex]
\left(\|S_{n+2}u_0-S_{n+1}u_0\|_{H^{s-1}}+\sum_{j=1}^5(\frac{1+\nu t}{\nu})^{1-\frac{1}{\rho_j}}
\|F_j\|_{\tilde{L}^{\rho_j}_t(H^{s-3+\frac{2}{\rho_j}})}\right)\\[1ex]
\leq 2C_0\left(\|\Delta_{n+1}u_0\|_{H^{s-1}}+
2\nu^{-1}\sum_{j=1}^5\|F_j\|_{\tilde{L}^{\rho_j}_t(H^{s-3+\frac{2}{\rho_j}})}\right).
\end{array}
\end{align}
Taking $\rho_1=2$, by Lemmas \ref{l11}-\ref{l12}, we have

\begin{align}
\begin{array}{l}
\|F_1\|_{\tilde{L}^2_t(H^{s-2})}=\|(u_n-u_{n-1})\cdot\nabla u_n\|_{\tilde{L}^2_t(H^{s-2})}\\[1ex]
\leq C\|T_{u_n-u_{n-1}}\nabla u_n+T_{\nabla u_n}(u_n-u_{n-1})+R(\nabla u_n,u_n-u_{n-1})\|_{\tilde{L}^2_t(H^{s-2})}\\[1ex]
\leq C\|u_n-u_{n-1}\|_{L^2_t(B^{s-2}_{\infty,\infty})}\|\nabla u_n\|_{L^\infty_t(H^0)}+
\|\nabla u_n\|_{L^\infty_t(B^{-1}_{\infty,\infty})}\|u_n-u_{n-1}\|_{L^2_t(H^{s-1})}\\[1ex]
+C\|\nabla u_n\|_{L^\infty_t(H^0)}\|u_n-u_{n-1}\|_{L^2_t(H^{s-1})}\\[1ex]
\leq CT_1^{\frac{1}{2}}\|u_n\|_{L^\infty_T(H^s)}\|u_n-u_{n-1}\|_{L^\infty_t(H^{s-1})}\\[1ex]
\leq C T_1^{\frac{1}{2}}E_1\|u_n-u_{n-1}\|_{L^\infty_t(H^{s-1})}.
\end{array}
\end{align}
Taking $\rho_2=2$ as well, we have
\begin{align}
\begin{array}{l}
\|F_2\|_{\tilde{L}^2_t(H^{s-2})}=\|\nabla(h_n-h_{n-1})\|_{\tilde{L}^2_t(H^{s-2})}\leq
CT_1^{\frac{1}{2}}\|h_n-h_{n-1}\|_{\tilde{L}^{\infty}_t(H^{s-1})}.
\end{array}
\end{align}
Taking $\rho_3=\frac{2}{s}$, by Lemmas \ref{l11}-\ref{l12}, we get
\begin{align}
\begin{array}{l}
\|F_3\|_{\tilde{L}^{\frac{2}{s}}_t(H^{s-3+s})}=\|\nu\nabla(\ln(1+h_n))\nabla(u_n-u_{n-1})\|_{L^{\frac{2}{s}}_t(H^{2s-3})}\\[1ex]
\leq C\|T_{\nabla(\ln(1+h_n))}\nabla(u_n-u_{n-1})\|_{L^{\frac{2}{s}}_t(H^{2s-3})}+C\|T_{\nabla u_n-u_{n-1}}\nabla(\ln(1+h_n))\|_{L^{\frac{2}{s}}_t(H^{2s-3})}\\[1ex]
+C\|R(\nabla(\ln(1+h_n)),\nabla(u_n-u_{n-1)})\|_{L^{\frac{2}{s}}_t(H^{2s-3})}\\[1ex]
\leq C\|\nabla(\ln(1+h_n))\|_{L^\infty_t(B^{s-2}_{\infty,\infty})}\|\nabla(u_n-u_{n-1})\|_{L^{\frac{2}{s}}_t(H^{s-1})}\\
+C\|\nabla(u_n-u_{n-1})\|_{L^{\frac{2}{s}}_t(B^{s-2}_{\infty,\infty})}\|\nabla(\ln(1+h_n))\|_{L^\infty_t(H^{s-1})}\\
+C\|\nabla(u_n-u_{n-1})\|_{L^{\frac{2}{s}}_t(H^{s-1})}\|\nabla(\ln(1+h_n))\|_{L^\infty_t(H^{s-1})}\\[1ex]
\leq C\|\ln(1+h_n)\|_{L^\infty_t(H^s)}\|u_n-u_{n-1}\|_{L^{\frac{2}{s}}_t(H^s)}\\[1ex]
\leq CT_1^{\frac{s-1}{2}}\|h_n\|_{L^\infty_t(H^s)}\|u_n-u_{n-1}\|_{L^2_t(H^s)}\\[1ex]
\leq CT_1^{\frac{s-1}{2}}E_2\|u_n-u_{n-1}\|_{L^2_t(H^s)}.
\end{array}
\end{align}

Then we deal with $F_4$ by the similar method. Taking $\rho_4=2$, we have
\begin{align}
\begin{array}{l}
\|F_4\|_{L^2_t(H^{s-2})}=\|\nu\frac{\nabla(h_n-h_{n-1})}{1+h_n}\nabla u_{n-1}\|_{L^2_t(H^{s-2})}\\[1ex]
\leq C\|\nabla u_{n-1}\nabla(h_n-h_{n-1})\|_{L^2_t(H^{s-2})}+C\|\frac{h_n}{1+h_n}\nabla
u_{n-1}\nabla(h_n-h_{n-1})\|_{L^2_t(H^{s-2})}\\[1ex]
=F_{41}+F_{42}.
\end{array}
\end{align}
By Lemmas \ref{l11}-\ref{l12}, we have
\begin{align}
\begin{array}{l}
F_{41}\leq C\|\nabla(h_n-h_{n-1})\|_{L^\infty_t(B^{s-3}_{\infty,\infty})}\|\nabla u_{n-1}\|_{L^2_T(H^1)}+
\|\nabla u_{n-1}\|_{L^2_t(L^\infty)}\|\nabla(h_n-h_{n-1})\|_{L^\infty_t(H^{s-2})}\\[1ex]
+C\|\nabla u_{n-1}\|_{L^2_t({H^1})}\|\nabla(h_n-h_{n-1})\|_{L^\infty_t(H^{s-2})}\\[1ex]
\leq C\|\nabla u_n\|_{L^2_t(L^\infty\cap H^1)}\|h_n-h_{n-1}\|_{L^\infty_t(H^{s-1})}\\[1ex]
\leq C\|u_n-u_{n-1}\|_{L^2_t(H^\frac{s+3}{2})}\|h_n-h_{n-1}\|_{L^\infty_t(H^{s-1})}\\[1ex]
\leq CT_1^{\frac{s-1}{4}}\|u_{n-1}\|_{L^{\frac{4}{3-s}}_t(H^{\frac{s+3}{2}})}\|h_n-h_{n-1}\|_{L^\infty_t(H^{s-1})}\\[1ex]
\leq CT_1^{\frac{s-1}{4}}E_1\|h_n-h_{n-1}\|_{L^\infty_t(H^{s-1})}.
\end{array}
\end{align}
Similarly as $F_{41}$, using Lemmas \ref{l11}-\ref{l12} again, we
have
\begin{align}
\begin{array}{l}
F_{42}\leq C\|\frac{h_n}{1+h_n}\nabla u_{n-1}\|_{L^2_t(L^\infty\cap H^1)}\|h_n-h_{n-1}\|_{L^\infty_t(H^{s-1})}\\[1ex]
\leq C\|\frac{h_n}{1+h_n}\nabla u_{n-1}\|_{L^2_t(H^{\frac{s+1}{2}})}\|h_n-h_{n-1}\|_{L^\infty_t(H^{s-1})}\\[1ex]
\leq C\|\frac{h_n}{1+h_n}\|_{L^\infty_t(H^s)}\|\nabla u_{n-1}\|_{L^2_t(H^{\frac{s+1}{2}})}\|h_n-h_{n-1}\|_{L^\infty_t(H^{s-1})}\\[1ex]
\leq CT_1^{\frac{s-1}{4}}\|h_n\|_{L^\infty_t(H^s)}\|u_{n-1}\|_{L^{\frac{4}{3-s}}_t(H^{\frac{s+3}{2}})}\|h_n-h_{n-1}\|_{L^\infty_t(H^{s-1})}\\[1ex]
\leq CT_1^{\frac{s-1}{4}}E_1E_2\|h_n-h_{n-1}\|_{L^\infty_t(H^{s-1})}.
\end{array}
\end{align}

For $F_5$, taking $\rho_5=2$, then we have
\begin{align}
\begin{array}{l}
\|F_5\|_{L^2_t(H^{s-2})}=\|\nu
(\frac{1}{1+h_n}-\frac{1}{1+h_{n-1}})\nabla h_{n-1}\nabla u_{n-1}\|_{L^2_t(H^{s-2})}\\[1ex]
=\nu\|(h_{n-1}-h_n)\frac{1}{1+h_n}\nabla(\ln(1+h_{n-1}))\nabla u_{n-1}\|_{L^2_t(H^{s-2})}\\[1ex]
\leq C\|(h_{n-1}-h_n)\nabla(\ln(1+h_{n-1}))\nabla u_{n-1}\|_{L^2_t(H^{s-2})}\\[1ex]
+\|(h_{n-1}-h_n)\frac{h_n}{1+h_n}\nabla(\ln(1+h_{n-1}))\nabla u_{n-1}\|_{L^2_t(H^{s-2})}\\[1ex]
=F_{51}+F_{52}.
\end{array}
\end{align}
By Lemmas \ref{l11}-\ref{l12}, we have
\begin{align}
\begin{array}{l}
F_{51}\leq C\|T_{\nabla(\ln(1+h_{n-1}))\nabla u_{n-1}}(h_n-h_{n-1})\|_{L^2_t(H^{s-2})}\\[1ex]
+C\|T_{h_n-h_{n-1}}(\nabla(\ln(1+h_{n-1}))\nabla u_{n-1})\|_{L^2_t(H^{s-2})}\\[1ex]
+C\|R(\nabla(\ln(1+h_{n-1}))\nabla u_{n-1},(h_n-h_{n-1}))\|_{L^2_t(H^{s-2})}\\[1ex]
\leq C\|\nabla(\ln(1+h_{n-1}))\nabla u_{n-1}\|_{L^2_t(B^{-1}_{\infty,\infty})}\|h_n-h_{n-1}\|_{L^\infty_t(H^{s-1})}\\[1ex]
+ C\|h_n-h_{n-1}\|_{L^\infty_t(B^{s-2}_{\infty,\infty})}\|\nabla(\ln(1+h_{n-1}))\nabla u_{n-1}\|_{L^2_t(H^0)}\\[1ex]
+\|\nabla(\ln(1+h_{n-1}))\nabla u_{n-1}\|_{L^2_t(H^0)}\|h_n-h_{n-1}\|_{L^\infty_t(H^{s-1})}\\[1ex]
\leq C\|\nabla(\ln(1+h_{n-1}))\nabla u_{n-1}\|_{L^2_t(H^0)}\|h_n-h_{n-1}\|_{L^\infty_t(H^{s-1})}.
\end{array}
\end{align}
By Lemmas \ref{l11}-\ref{l12} again, we have
\begin{align}
\begin{array}{l}
\|\nabla(\ln(1+h_{n-1}))\nabla u_{n-1}\|_{L^2_t(H^0)}\\[1ex]
\leq C\|\nabla(\ln(1+h_{n-1}))\|_{L^\infty_t(B^{s-2}_{\infty,\infty})}\|\nabla u_{n-1}\|_{L^2_t(H^{2-s})}\\[1ex]
+C\|\nabla u_{n-1}\|_{L^2_t(B^{1-s}_{\infty,\infty})}\|\nabla(\ln(1+h_{n-1}))\|_{L^\infty_t(H^{s-1})}\\[1ex]
+C\|\nabla(\ln(1+h_{n-1}))\|_{L^\infty_t(H^{s-1})}\|\nabla u_{n-1}\|_{L^2_t(H^{2-s})}\\[1ex]
\leq C\|h_{n-1}\|_{L^\infty_t(H^s)}\|u_{n-1}\|_{L^2_t(H^{3-s})}\\[1ex]
\leq CT_1^{\frac{s-1}{4}}\|h_{n-1}\|_{L^\infty_t(H^s)}\|u_{n-1}\|_{L^{\frac{4}{3-s}}_t(H^{\frac{s+3}{2}})}.
\end{array}
\end{align}
Combining (3.22) and (3.23), we have
\begin{align}
\begin{array}{l}
F_{51}\leq
CT_1^{\frac{s-1}{4}}\|h_{n-1}\|_{L^\infty_t(H^s)}\|u_{n-1}\|_{L^{\frac{4}{3-s}}_t(H^{\frac{s+3}{2}})}
\|h_n-h_{n-1}\|_{L^\infty_t(H^{s-1})}\\[1ex]
\leq CT_1^{\frac{s-1}{4}}E_1E_2\|h_n-h_{n-1}\|_{L^\infty_t(H^{s-1})}.
\end{array}
\end{align}
Similarly, we have
\begin{align}
\begin{array}{l}
F_{52}\leq C\|\frac{h_n}{1+h_n}\nabla(\ln(1+h_{n-1}))\nabla u_{n-1}\|_{L^2_t(H^0)}\|h_n-h_{n-1}\|_{L^\infty_t(H^{s-1})}\\[1ex]
\leq C\|\frac{h_n}{1+h_n}\|_{L^\infty_t(H^s)}\|\nabla(\ln(1+h_{n-1}))\nabla u_{n-1}\|_{L^2_t(H^0)}\|h_n-h_{n-1}\|_{L^\infty_t(H^{s-1})}\\[1ex]
\leq CT_1^{\frac{s-1}{4}}\|h_{n-1}\|^2_{L^\infty_t(H^s)}\|u_{n-1}\|_{L^{\frac{4}{3-s}}_t(H^{\frac{s+3}{2}})}
\|h_n-h_{n-1}\|_{L^\infty_t(H^{s-1})}\\[1ex]
\leq CT_1^{\frac{s-1}{4}}E_1E_2^2\|h_n-h_{n-1}\|_{L^\infty_t(H^{s-1})}.
\end{array}
\end{align}
Thus, we have
\begin{align}
\begin{array}{l}
\|u_{n+1}-u_n\|_{\tilde{L}^\infty_t(H^{s-1})}
\leq 2C_0\|\Delta_{n+1}u_0\|_{H^{s-1}}+CT^{\frac{1}{2}}E_1\|u_n-u_{n-1}\|_{L^\infty_t(H^{s-1})}\\[1ex]
+CT_1^{\frac{1}{2}}\|h_n-h_{n-1}\|_{\tilde{L}^{\infty}_t(H^{s-1})}
+CT_1^{\frac{s-1}{2}}E_2\|u_n-u_{n-1}\|_{L^2_t(H^s)}\\[1ex]
+CT_1^{\frac{s-1}{4}}E_1\|h_n-h_{n-1}\|_{L^\infty_t(H^{s-1})}
+CT_1^{\frac{s-1}{4}}E_1E_2\|h_n-h_{n-1}\|_{L^\infty_t(H^{s-1})}\\[1ex]
+CT_1^{\frac{s-1}{4}}E_1E_2\|h_n-h_{n-1}\|_{L^\infty_t(H^{s-1})}+CT^{\frac{s-1}{4}}E_1E_2^2\|h_n-h_{n-1}\|_{L^\infty_t(H^{s-1})}\\[1ex]
\leq 2C_0\|\Delta_{n+1}u_0\|_{H^{s-1}}+CT_1^{\frac{s-1}{4}}E_1\|u_n-u_{n-1}\|_{L^\infty_t(H^{s-1})}+CT_1^{\frac{s-1}{4}}E_2\|u_n-u_{n-1}\|_{L^2_t(H^s)}\\[1ex]
+CT_1^{\frac{s-1}{4}}(1+E_1+E_1E_2+E_1E_2^2)\|h_n-h_{n-1}\|_{L^\infty_t(H^{s-1})}.
\end{array}
\end{align}

For $h_{n+1}-h_n$, we have
\begin{align}
\begin{array}{l}
\|h_{n+1}-h_n\|_{\tilde{L}^{\infty}_t(H^{s-1})}\leq e^{C_0\int_0^t\|\nabla u_n\|_{H^1\cap L^\infty}}\times
\big(\|\Delta_{n+1}h_0\|_{H^{s-1}}
+\|\sum_{j=1}^4G_j\|_{\tilde{L}^1_t(H^{s-1})}\big)\\[1ex]
\leq2\big(\|\Delta_{n+1}h_0\|_{H^{s-1}}
+T_1^{\frac{1}{2}}\|\sum_{j=1}^4J_j\|_{\tilde{L}^2_t(H^{s-1})}\big).
\end{array}
\end{align}

From Lemma \ref{l12}, we have
\begin{align}
\begin{array}{l}
\|J_1\|_{\tilde{L}_t^2(H^{s-1})}=\|(u_n-u_{n-1})\cdot\nabla h_n\|_{\tilde{L}_t^2(H^{s-1})}\\[1ex]
\leq C\|u_n-u_{n-1}\|_{L^2_t(L^\infty)}\|\nabla h_n\|_{L^\infty_t(H^{s-1})}+C\|\nabla h_n\|_{L^\infty_t(B^{s-2}_{\infty,\infty})}\|u_n-u_{n-1}\|_{L^2_t(H^1)}\\[1ex]
+C\|\nabla h_n\|_{L^\infty_t(H^{s-1})}\|u_n-u_{n-1}\|_{L^2_t(H^1)}\\[1ex]
\leq C\|h_n\|_{L^\infty_t(H^s)}\|u_n-u_{n-1}\|_{L^2_t(H^s)}\\[1ex]
\leq CE_2\|u_n-u_{n-1}\|_{L^2_t(H^s)}.
\end{array}
\end{align}
Clearly, we have
\begin{align}
\begin{array}{l}
\|J_2\|_{\tilde{L}_t^2(H^{s-1})}=\|div(u_n-u_{n-1})\|_{L^2_t(H^{s-1})}\leq C\|u_n-u_{n-1}\|_{L^2_t(H^s)}.
\end{array}
\end{align}

By Lemma \ref{l12}, we obtain
\begin{align}
\begin{array}{l}
\|J_3\|_{\tilde{L}_t^2(H^{s-1})}=\|h_n\,div\,(u_n-u_{n-1})\|_{\tilde{L}_t^2(H^{s-1})}\\[1ex]
\leq C\|h_n\|_{L^\infty_t(L^\infty)}\|div\,(u_n-u_{n-1})\|_{\tilde{L}_t^2(H^{s-1})}
+C\|div\,(u_n-u_{n-1})\|_{\tilde{L}_t^2(B^{s-2}_{\infty,\infty})}\|h_n\|_{L^\infty_t(H^1)}\\[1ex]
+C\|h_n\|_{L^\infty_t(H^1)}\|div\,(u_n-u_{n-1})\|_{\tilde{L}_t^2(H^{s-1})}\\[1ex]
\leq C\|h_n\|_{L^\infty_t(H^s)}\|u_n-u_{n-1}\|_{\tilde{L}_t^2(H^s)}\\[1ex]
\leq CE_2\|u_n-u_{n-1}\|_{\tilde{L}_t^2(H^s)}.
\end{array}
\end{align}
Similarly to $J_3$, by Lemma \ref{l12}, we have
\begin{align}
\begin{array}{l}
\|J_4\|_{\tilde{L}_t^2(H^{s-1})}=\|(h_n-h_{n-1})\,div\,u_{n-1}\|_{\tilde{L}_t^2(H^{s-1})}\\[1ex]
\leq C\|h_n-h_{n-1}\|_{L^\infty_t(B^{s-2}_{\infty,\infty})}\|div\,u_{n-1}\|_{L^2_t(H^1)}+
\|div\,u_{n-1}\|_{L^2_t(L^\infty)}\|h_n-h_{n-1}\|_{L^\infty_t(H^{s-1})}\\[1ex]
+C\|div\,u_{n-1}\|_{L^2_t(H^1)}\|h_n-h_{n-1}\|_{L^\infty_t(H^{s-1})}\\[1ex]
\leq C\|u_{n-1}\|_{L^2_t(H^{s+1})}\|h_n-h_{n-1}\|_{L^\infty_t(H^{s-1})}\\[1ex]
\leq CE_1\|h_n-h_{n-1}\|_{L^\infty_t(H^{s-1})}.
\end{array}
\end{align}
Thus we have
\begin{align}
\begin{array}{l}
\|h_{n+1}-h_n\|_{\tilde{L}^{\infty}_t(H^{s-1})}
\leq2\big(\|\Delta_{n+1}h_0\|_{H^{s-1}}
+T_1^{\frac{1}{2}}\|\sum_{j=1}^4G_j\|_{\tilde{L}^2_t(H^{s-1})}\big)\\[1ex]
\leq 2\|\Delta_{n+1}h_0\|_{H^{s-1}}+CT_1^{\frac{1}{2}}(1+E_2)\|u_n-u_{n-1}\|_{L^2_t(H^s)}+CT_1^{\frac{1}{2}}E_1
\|h_n-h_{n-1}\|_{L^\infty_t(H^{s-1})}.
\end{array}
\end{align}
Similarly as $\|u_{n+1}-u_n\|_{\tilde{L}^\infty_t(H^{s-1})}$, for all $\rho_j\leq2$, we have
\begin{align}
\begin{array}{l}
\|u_{n+1}-u_n\|_{L^2_t(H^s)}\leq \nu^{-\frac{1}{2}}C_0e^{C_0\int_0^t\|u_n\|_{H^{s+1}}}dt'\Big((1+\nu t)^{\frac{1}{2}}\|\Delta_{n+1}u_0\|_{H^{s-1}}\\[1ex]
+\sum_{j=1}^5(1+\nu t)^{1+\frac{1}{2}-\frac{1}{\rho_j}}\nu^{\frac{1}{\rho_j}-1}\|F_j\|_{\tilde{L}^{\rho_j}_t(H^{s-3+\frac{2}{\rho_j}})}\Big)\\[1ex]
\leq 4\nu^{-1}C_0\|\Delta_{n+1}u_0\|_{H^{s-1}}+4\nu^{-1}\sum_{j=1}^5\|F_j\|_{\tilde{L}^{\rho_j}_t(H^{s-3+\frac{2}{\rho_j}})}\\[1ex]
\leq4\nu^{-1}C_0\|\Delta_{n+1}u_0\|_{H^{s-1}}+CT_1^{\frac{s-1}{4}}E_1\|u_n-u_{n-1}\|_{L^\infty_t(H^{s-1})}+CT_1^{\frac{s-1}{4}}E_2\|u_n-u_{n-1}\|_{L^2_t(H^s)}\\[1ex]
+CT_1^{\frac{s-1}{4}}(1+E_1+E_1E_2+E_1E_2^2)\|h_n-h_{n-1}\|_{L^\infty_t(H^{s-1})}.
\end{array}
\end{align}
We also have
\begin{align}
\begin{array}{l}
\|\Delta_{n+1}u_0\|_{H^{s-1}}\leq2^{-(n+1)}\|\Delta_{n+1}u_0\|_{H^s}
\leq2^{-(n+1)}\|u_0\|_{H^s}\leq 2^{-(n+1)}E_1,\\[1ex]
\|\Delta_{n+1}h_0\|_{H^{s-1}}\leq2^{-(n+1)}\|\Delta_{n+1}h_0\|_{H^s}
\leq2^{-(n+1)}\|h_0\|_{H^s}\leq 2^{-(n+1)}E_2.
\end{array}
\end{align}
We will complete this proof by induction.
By Proposition \ref{p1}, we have
$$\|u_2-u_1\|_{\tilde{L}_{T_1}^\infty(H^{s-1})}\leq\|u_1\|_{\tilde{L}_{T_1}^\infty(H^{s-1})}+
\|u_2\|_{\tilde{L}_{T_1}^\infty(H^{s-1})}$$
$$\leq\|u_1\|_{\tilde{L}_{T_1}^\infty(H^s)}+
\|u_2\|_{\tilde{L}_{T_1}^\infty(H^s)}$$
$$\leq 2E_1.$$
Similarly, we have
$$\|u_2-u_1\|_{\tilde{L}_{T_1}^2(H^s)}\leq 2E_1,$$
$$\|h_2-h_1\|_{\tilde{L}_{T_1}^\infty(H^{s-1})}\leq 2E_2.$$
We assume that
$$\|u_n-u_{n-1}\|_{\tilde{L}_{T_1}^\infty(H^{s-1})} \leq8\times2^{-n}E_1,$$
$$\|u_n-u_{n-1}\|_{\tilde{L}_{T_1}^2(H^s)}\leq8\times 2^{-n}E_1,$$
$$\|h_n-h_{n-1}\|_{\tilde{L}_{T_1}^\infty(H^{s-1})}\leq8\times 2^{-n}E_2.$$
Thus take a suitable $T_1$ satisfying
\begin{align}
\left\{
\begin{array}{l}
8CT_1^{\frac{s-1}{4}}E_1\leq1,~8CT_1^{\frac{s-1}{4}}E_2\leq1,\\
8CT_1^{\frac{s-1}{4}}(1+E_1+E_1E_2+E_1E_2^2)E_2E_1^{-1}\leq1,\\
8CT_1^{\frac{1}{2}}E_1\leq1,~8CT_1^{\frac{1}{2}}E_1(1+E_2)E_2^{-1}\leq1.
\end{array}
\right.
\end{align}
Where $C$ in subsection 3.2 only depends on $C_0,\nu$ as well.
Thus using (3.26), (3.32), (3.33), we can obtain that
$$\|u_{n+1}-u_n\|_{\tilde{L}_{T_1}^\infty(H^{s-1})} \leq8\times2^{-(n+1)}E_1,$$
$$\|u_{n+1}-u_n\|_{\tilde{L}_{T_1}^2(H^s)}\leq8\times 2^{-(n+1)}E_1,$$
$$\|h_{n+1}-h_n\|_{\tilde{L}_{T_1}^\infty(H^{s-1})}\leq8\times 2^{-(n+1)}E_2.$$
This completes the proof of Proposition \ref{p2}.
\subsection{Existence and uniqueness of local solutions}
~~~~~In this subsection, we investigate the uniqueness of the local solution to the system (2.1). By Proposition \ref{p2}, the approximative sequence $(u_n,h_n)$ of the problem (3.1) is a Cauchy sequence in $\chi([0,T_1],s-1,E_1,E_2)$ with $s\in(1,2)$. So the limit (u,h) is a solution of the Cauchy problem (2.1). From Proposition \ref{p1}, we obtain that this sequence is bounded in $\chi([0,T],s,E_1,E_2)$, so it's also the Cauchy sequence in $\chi([0,T_1],s',E_1,E_2)$ for all $s'<s$ by interpolation, and the limit is in $\chi([0,T_1],s,E_1,E_2)$. Thus we have proved local existence result in Theorem \ref{t1}.

For the uniqueness result in Theorem \ref{t1}, let $(u,h)$ and
$(v,g)$ satisfy the problem (2.1) with the initial data
$(u_0,h_0),~(v_0,h_0)\in H^s\times H^s$ respectively. Then we have
\begin{align}
\left\{
\begin{array}{l}
\partial_t(u-v)+u\cdot\nabla(u-v)-\nu\Delta(u-v)=G_1(u,h)-G_1(v,g),\\
\partial_t(h-g)+u\cdot\nabla(h-g)=(u-v)\nabla g+G_2(u,h)-G_2(v,g),\\
(u-v)|_{t=0}=0, (h-g)|_{t=0}=0.
\end{array}
\right.
\end{align}
Using Lemmas \ref{l14}-\ref{l15}, similarly as the proof of Proposition \ref{p2}, for all $t\leq T_1$, we can get,
\begin{align}
\begin{array}{l}
\|u-v\|_{\tilde{L}^{\infty}(H^s)}+\|u-v\|_{\tilde{L}^2(H^{s+1})}+\|h-g\|_{\tilde{L}^{\infty}(H^s)}\\
\leq Ce^{Ct}(\|u_0-v_0\|_{H^s}+\|h_0-g_0\|_{H^s})\\[1ex]
+Ct^{\alpha_1}\|u-v\|_{\tilde{L}^{\infty}(H^s)}+Ct^{\alpha_2}\|u-v\|_{\tilde{L}^2(H^{s+1})}+
Ct^{\alpha_3}\|h-g\|_{\tilde{L}^{\infty}(H^s)},
\end{array}
\end{align}
where ${\alpha_1},{\alpha_2},{\alpha_3}\in(0,1)$, and $C$ only depends on $C_0,~\nu,~E_1,~E_2$.
Taking a suitable $T_2\leq T_1$ gives the uniqueness of Theorem \ref{t1} in $[0,T_2]$, and thus we get the uniqueness in $[0,T_1]$.

\section{Global existence}

In \cite{W}, Wang and Xu got the global existence of the Cauchy
problem (1.1) for small initial data in $H^s\times H^s$ with
$s>2$, since they showed only the existence of the local solutions
in $H^s\times H^s$ with $s>2$. And They also pointed out that (see
Remark of Theorem 4.2 in \cite{W}): if one could establish the
local solutions in $H^s\times H^s$ with $s>1$, one can also obtain
the global solution with initial data in $H^s\times H^s$ with
$s>1$. Since the proof of global existence is the same as that of
the $s>2$ case \cite{W}, here we omit it.
\newline
\textbf{Remark 4.1} Note Wang and Xu in \cite{W} established the
local well-posedness of the system (1.1) for any initial data and
got the global solutions to the system (1.1) for small initial
data in Sobolev spaces $H^s$, $s>2$. Our obtained result in
Theorem 1.1 shows that Wang and Xu's conjecture in \cite{W} holds
true.

\vspace*{2em} \noindent\textbf{Acknowledgements}.This work was
partially supported by NNSFC (No.11271382 No. 10971235), RFDP (No.
20120171110014), and the key project of Sun Yat-sen University.

\phantomsection
\addcontentsline{toc}{section}{\refname}

\end{document}